\documentclass[11pt]{article}
\usepackage{}
\usepackage{graphicx}
\usepackage{graphicx,subfigure}
\usepackage{epsfig}
\usepackage{amssymb,pict2e,color}
\usepackage{mathrsfs}
\usepackage{amsmath}
\usepackage{ifpdf}
\usepackage [latin1]{inputenc}

\newtheorem{theorem}{Theorem}[section]

\newtheorem{lemma}[theorem]{Lemma}

\newtheorem{claim}[theorem]{Claim}
\newtheorem{proposition}[theorem]{Proposition}

\def\whitebox{{\hbox{\hskip 1pt
 \vrule height 6pt depth 1.5pt
 \lower 1.5pt\vbox to 7.5pt{\hrule width
    3.2pt\vfill\hrule width 3.2pt}%
 \vrule height 6pt depth 1.5pt
 \hskip 1pt } }}
\def\qed{\ifhmode\allowbreak\else\nobreak\fi\hfill\quad\nobreak
     \whitebox\medbreak}
\newcommand{\proof}{\noindent{\it Proof.}\ }

\newcommand{\ignore}[1]{}

\begin {document}

\baselineskip 16pt
\title{The multiplicative Zagreb indices of graphs with  given connectivity or number of pendant vertices}

 \author{\small Shengjin Ji$^{\dagger,\ddagger,} \footnote{ Corresponding author. \newline Email addresses:    jishengjin2013@163.com(S. Ji), shaohuiwang@yahoo.com(S. Wang), muchet@savannahstate.edu(T. Muche), sakander@mail.ustc.edu.cn(S. Hayat).}$, Shaohui Wang$^{\nmid}$, Tilahun Muche$^\nmid$, Sakander Hayat$^\nparallel$\\
 \small  $^{\dagger}$ School of Mathematics and Statistics, Shandong University of Technology, Zibo, Shandong 255049, China\\
 \small  $^{\nmid}$ Department of Mathematics, Savannah State University, Savannah, GA 31404, USA\\
 \small $^{\ddagger}$ School of Mathematics, Shandong University,
Jinan, Shandong 250100, China\\
\small $^{\nparallel}$ School of Mathematical Sciences, University of Science and Technology of China, Hefei, 230026, China
}

\date{}
\maketitle
\begin{abstract}

For a graph $G$, the first multiplicative Zagreb index $\prod_1(G) $ is the product of squares of vertex degrees, and the second multiplicative Zagreb index $\prod_2(G) $ is the product of  products of degrees of pairs of adjacent vertices.   In this paper, we explore graphs with extremal $\Pi_{1}(G)$ and $\Pi_{2}(G)$ in terms of (edge) connectivity and pendant vertices. The corresponding extremal graphs are characterized with given connectivity at most $k$ and $p$ pendant vertices.  In addition,  the maximum and minimum values of $\prod_1(G) $ and $\prod_2(G) $  are provided. Our results extend and enrich some known conclusions.

\end{abstract}

\section{Introduction}
A topological index is a single number which can be used to describe some properties of  a molecular graph that is  a finite simple graph, representing the carbon-atom skeleton
of an organic molecule of a hydrocarbon.
In recent decades, these numerical quantities have been found   useful for the study of quantitative structure-property relationships (QSPR) and quantitative
structure-activity relationships (QSAR)  and for  the structural essence of biological and chemical compounds.  The well-known Randi\'{c} index is one of the most important topological indices.

In 1975, Randi\'{c} introduced a  moleculor quantity of branching index \cite{1}, which has been known as the famous Randi\'{c} connectivity index and that is a most useful structural descriptor in QSPR and QSAR, see \cite{2,3,4,5}.  Mathematicians have considerable interests in the structural and applied aspects of  Randi\'{c} connectivity index, see \cite{6,7,8,9}. Based on the successful considerations,  Zagreb indices \cite{10} are introduced as an expected formula for the total $\pi$-electron energy of conjugated molecules as follows.
\begin{eqnarray} \nonumber
M_1(G) = \sum_{u \in V(G)} d(u)^2
~\text{ and}
~ M_2(G) = \sum_{uv \in E(G)} d(u)d(v),
\end{eqnarray}
where $G$ is a (molecular) graph, $uv$ is a bond between two atoms $u$ and $v$, and $d(u)$ (resp. $d(v)$) is the number of atoms that are connected with $u$ (resp. $v$).
Zagreb indices have also been employed as molecular descriptors in QSPR and QSAR, see \cite{11,12}.  Recently, Todeschini et al. (2010) \cite{13,14} proposed the
following multiplicative variants of molecular structure descriptors:
\begin{eqnarray} \nonumber
 \prod_1(G) = \prod_{u \in V(G)} d(u)^2 ~
\text{ and}\;\;
 \prod_2(G) = \prod_{uv \in E(G)} d(u)d(v) = \prod_{u \in V(G)}
d(u)^{d(u)}.
\end{eqnarray}
In the interplay among mathemactics, chemistry and physics, it is not surprising that  there are numerous studies of properties of the (multiplicative) Zagreb indices of molecular graphs \cite{1000,1001,1002,1003,s1,s2,s3,1006,sx2015,mr2018,vb2020}.

In view of these results, researchers are interested in
finding upper and lower bounds for multiplicative Zagreb indices of graphs  and characterizing the graphs
in which the maximal and minimal index values are attained. In view of the   above problems,  various mathematical and computational properties of Zagreb
indices have been investigated in \cite{16, 17, 18}. Other directions of investigation include studies of relation between multiplicative Zagreb indices and the corresponding invariant of elements of the graph $G$ (vertices, pendant vertices, diameter, maximum degree, girth, cut edge, cut vertex, connectivity, perfect matching).

For instance, the first and second multiplicative Zagreb indices for a class of chemical dendrimers are explored by Iranmanesh et al.~\cite{Iranmanesh20102}.
 Considering trees,  unicyclic
graphs and bicyclic graphs, Borovi\'canin et al.~\cite{Borov2016} introduced the
bounds on Zagreb indices with a fixed  domination number.
 The maximum and minimum
Zagreb indices of trees with given
number of vertices of maximum degree are proposed by  Borovi$\acute{c}$anin and Lampert \cite{Bojana2015}.  Xu and Hua~\cite{Xu20102} introduced a unified approach to
characterize  maximal and minimal multiplicative Zagreb
indices, respectively.
Considering the trees of higher dimension, i.e. $k$-trees, Wang and Wei~\cite{Wang2015} provided the maximum and minimum values of these indices and   the corresponding extremal graphs.
 Some sharp upper bounds for
$\prod_1$-index and $\prod_2$-index in terms of graph parameters are investigated by Liu and Zhang \cite{Liuz20102},
including  the order, size and radius of graphs.
Ji and Wang \cite{WW} provided the sharp lower bounds of Zagreb indices of graphs with given
number of cut vertices.
  The bounds for
the moments and the probability generating function of these
indices in a randomly chosen molecular graph with tree structure
of given order  are studied by Kazemi~\cite{Ramin2016}.  Li and Zhao obtained sharp upper bounds on Zagreb indices of bicyclic graphs with a given matching number~\cite{LL,zhaoli}.

In light of the information available for multiplicative Zagreb indices,
and inspired by  above results, in this paper we further
investigate these indices of graphs with  a given (edge) connectivity and number of pendant vertices. We give some basic
properties of the first and the second multiplicative Zagreb indices.  The maximum and minimum values of $\prod_1(G) $ and $\prod_2(G) $ of graphs with given (edge) connectivity at most $k$ and $p$ pendant vertices are provided.  In addition, the
corresponding extremal graphs are charaterized.
 In our exposition, we will use the terminology and notations of (chemical) graph theory (see \cite{BB,NT}).

\section{Preliminaries}
Let $G$ be a simple  connected  graph, denoted by $G = (V(G), E(G))$, in which $V = V (G)$ is  vertex
set and $E = E(G)$ is  edge set. For a vertex $v \in V(G)$,   the neighborhood of  $v $ is the set $N(v) = N_G(v) = \{w \in V(G), vw \in E(G)\}$, and $d_G(v)$ (or   $d(v)$)  is the degree of $v$ with $ d_{G}(v) = |N(v)|$. For $i\geq 0$, $n_i$ denoted the number of vertices of degree $i$.
For $S \subseteq V(G)$ and  $F \subseteq E(G)$,   we use $G[S]$ for the subgraph of $G$ induced by the vertex set $S$, $G - S$ for the subgraph induced by  $V(G) - S$ and $G - F$ for the subgraph of $G$ obtained by deleting $F$. If $G-S$ contains at least 2 components, then $S$ is said to be a vertex cut set of $G$. Similarly, if $G-F$ contains at least 2 components, then $F$ is called an edge cut set.

A graph $G$ is said to be $k$-connected with $k \geq 1$, if either $G$ is complete graph $K_{k+1}$, or it has at least $k+2$ vertices and contains no $(k-1)$-vertex cut.
The connectivity of $G$, denoted by $\kappa(G)$, is defined as the maximal value of $k$ for which a connected graph $G$ is $k$-connected. Similarly,
for $k \geq 1,$ a graph $G$ is called $k$-edge-connected if it has at least two vertices and does not contain a $(k-1)$-edge cut. The maximal value of $k$ for which a connected graph $G$ is $k$-edge-connected is said to be  the edge connectivity
of $G$, denoted by $\kappa'(G)$. According to the above definitions, the following proposition is obtained.

\begin{proposition} \label{p3} Let $G$ be a graph with $n$ vertices. Then \\
(i) $\kappa(G) \leq \kappa'(G) \leq n-1,$
\\
(ii) $\kappa(G)=\kappa'(G)=n-1$ if and only if $G\cong K_n$.
\end{proposition}
Let $\mathbb{V}_n^k$ be the set of graphs with $n$ vertices and $\kappa(G) \leq k \leq n-1.$ Denote $\mathbb{E}_n^k$ by the set of graphs with $n$ vertices and $\kappa'(G) \leq k \leq n-1.$  Note that if $|V(G)| = n$  and $|E(G)| = n-1$,  then $G$ is a tree. Let $P_n$ and $S_n$ be special trees: a path and a star of $n$ vertices.  The graph $K_n^k$ is obtained by joining $k$ vertices of $K_{n-1}$ to an isolated vertex, see Fig 1. Then $K_n^k \in \mathbb{E}_n^k \subset \mathbb{V}_n^k$.

\vspace{2mm}

\begin{center}
\begin{picture}(327,116)\linethickness{0.8pt}
\cbezier(168.5,59.5)(169.7,-6.1)(204.3,-6.1)(205.5,59.5)
\cbezier(168.5,59.5)(169.7,125.1)(204.3,125.1)(205.5,59.5)
\put(171.1,46.6){\footnotesize$K_j$}
\Line(187.5,30.5)(243,37)
\Line(187.5,30.5)(243,51.5)
\Line(187.5,30.5)(243,67.5)
\Line(187.5,30.5)(243,82)
\Line(187,46)(243,37)
\Line(187,46)(243,51.5)
\Line(187,46)(243,67.5)
\Line(187,46)(243,82)
\Line(186.5,61)(243,37)
\Line(186.5,61)(243,51.5)
\Line(186.5,61)(243,67.5)
\Line(186.5,61)(243,82)
\Line(186.5,76.5)(243,37)
\Line(186.5,76.5)(243,51.5)
\Line(186.5,76.5)(243,67.5)
\Line(186.5,76.5)(243,82)
\Line(186,91.5)(243,37)
\Line(186,91.5)(243,51.5)
\Line(186,91.5)(243,67.5)
\Line(186,91.5)(243,82)
\put(236.1,19.6){\footnotesize$H_k$}
\cbezier(229,56.8)(230,-5.2)(258,-5.2)(259,56.8)
\cbezier(229,56.8)(230,118.7)(258,118.7)(259,56.8)
\Line(243,37)(303.5,20.5)
\Line(243,37)(303,38)
\Line(243,37)(303,52)
\Line(243,37)(303,68)
\Line(243,37)(302.5,82)
\Line(243,37)(302,97)
\Line(243,51.5)(303.5,20.5)
\Line(243,67.5)(303.5,20.5)
\Line(243,51.5)(303,38)
\Line(243,51.5)(303,52)
\Line(243,51.5)(303,68)
\Line(243,51.5)(302.5,82)
\Line(243,51.5)(302,97)
\Line(243,67.5)(303,38)
\Line(243,67.5)(303,52)
\Line(243,67.5)(303,68)
\Line(243,67.5)(302.5,82)
\Line(243,67.5)(302,97)
\Line(243,82)(303.5,20.5)
\Line(243,82)(303,38)
\Line(243,82)(303,52)
\Line(243,82)(303,68)
\Line(243,82)(302.5,82)
\Line(243,82)(302,97)
\cbezier(287,58.5)(288.3,-17.7)(325.7,-17.7)(327,58.5)
\cbezier(287,58.5)(288.3,134.7)(325.7,134.7)(327,58.5)
\put(305,50.1){\footnotesize$K_{n-k-j}$}
\put(21.6,48.1){\footnotesize$K_{n-1}$}
\put(45,24){$u_k$}
\put(45,90){$u_1$}
\put(67.3,44){$u_3$}
\put(67,65){$u_2$}
\put(133.8,57){$u$}
\put(41.5,58){\circle{82}}
\Line(138,54)(56.1,96.3)
\Line(138,54)(81.3,67.8)
\Line(138,54)(80.2,44.3)
\Line(138,54)(56,19.6)
\Line(56.1,96.3)(56,19.6)
\put(60,27){\circle*{2}}
\put(62.5,29.5){\circle*{2}}
\put(65,32){\circle*{2}}
\put(67.5,34.5){\circle*{2}}
\put(70,37){\circle*{2}}
\put(72.5,39.5){\circle*{2}}
\put(187.5,30.5){\circle*{4}}
\put(243,37){\circle*{4}}
\put(243,51.5){\circle*{4}}
\put(243,67.5){\circle*{4}}
\put(243,82){\circle*{4}}
\put(187,46){\circle*{4}}
\put(186.5,61){\circle*{4}}
\put(186.5,76.5){\circle*{4}}
\put(186,91.5){\circle*{4}}
\put(303.5,20.5){\circle*{4}}
\put(303,38){\circle*{4}}
\put(303,52){\circle*{4}}
\put(303,68){\circle*{4}}
\put(302.5,82){\circle*{4}}
\put(302,97){\circle*{4}}
\put(138,54){\circle*{4}}
\put(56.1,96.3){\circle*{4}}
\put(81.3,67.8){\circle*{4}}
\put(80.2,44.3){\circle*{4}}
\put(56,19.6){\circle*{4}}
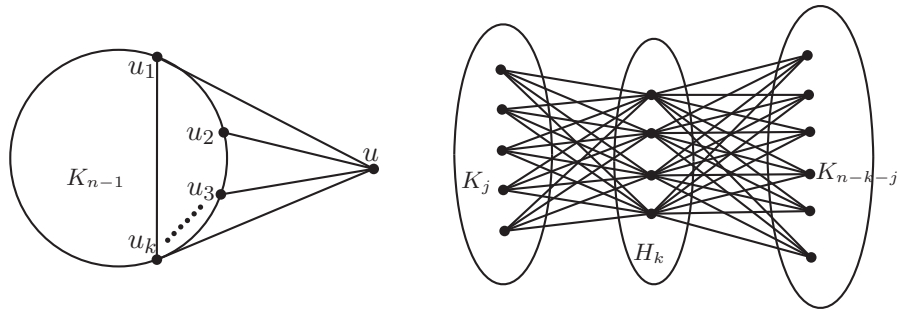
\put(160,-10){\makebox(0,0){Figure.1 The graphs $K_n^k$ and $G(j,n-k-j) = K_j \oplus H_k \oplus K_{n-k-j}$.}}
\end{picture}
\end{center}
\vspace{5mm}

Considering the concepts of $\prod_1(G)$ and $\prod_2(G)$, the  following proposition is routinely obtained.

\begin{proposition}\label{p4}
Let $e$ be an edge of a graph $G \in  \mathbb{V}_n^k$ (resp. $\mathbb{E}_{n}^{k}$). Then\\
(i) $G- e \in \mathbb{V}_n^k$ (resp. $\mathbb{E}_{n}^{k}$),\\
(ii) $\prod_i(G-e) < \prod_i({G}),  $i=1,2$.$
\end{proposition}

In addition, by elementary calculations, these three statements are deduced.

\begin{proposition}\label{p1}For $m\geq 0$,
  $F_1(x) = \frac{(x+m)^x}{(x-1+m)^{x-1}}$ is monotonically increasing in $(0,+\infty).$
\end{proposition}

\begin{proposition}\label{p2}
 If $m \geq 0$, then  $F_2(x) = \frac{x^x}{(x+m)^{x+m}}$ is a decreasing function  in interval $(0,+\infty ).$
\end{proposition}
\begin{proposition}\label{p3}
  $F_3(x) = (x)^2(n-x)^2$ is  monotonically increasing for $x\in[1,\lfloor\frac{n}{2}\rfloor].$
\end{proposition}

\section{Lemmas}

We first provide some lemmas, which are  important in proving our main results.

\begin{lemma}\label{l1}\cite{Iranmanesh20102}
Let $T$ be a tree on n vertices. If $T$ is not $P_n$ or $S_n$,  then both $\prod_1{(T)} > \prod_1(S_n)$ and $\prod_2(T) >\prod_2(P_n)$ holds.
\end{lemma}

Considering the definitions of $\prod_1(G)$ and $\prod_2(G)$, we have the following lemma.
\begin{lemma}\label{l2}
Let $u, v \in V(G)$ such that $uv \notin E(G)$. Then $$\prod_1(G+uv) > \prod_1(G)~\text{and}~\prod_2(G+uv) > \prod_2(G).$$
\end{lemma}

Given two graphs $G_1$ and $G_2$, if $V(G_1) \cap V(G_2) = \phi$, then the join graph $G_1 \oplus  G_2$ is a graph with vertex set $V(G_1) \cup V(G_2)$ and edge set $E(G_1) \cup E(G_2) \cup \{uv, u\in V(G_1), v \in V(G_2)\}$.

\begin{lemma}\label{l3} Let $G(j,n-k-j) = K_j \oplus H_k \oplus K_{n-k-j}$ be a graph with $n$ vertices, in which $K_j$ and $K_{n-k-j}$ are cliques, and $H_k$ is a graph with $k$ vertices, see Fig 2. If $k \geq 1$ and $2 \leq j \leq \frac{n-k}{2}$, then
$$\prod_1(G(j,n-k-j)) < \prod_1(G(1, n-k-1)).$$
\end{lemma}

\begin{proof}
We consider the graph from $G_1= G(j,n-k-j)$ to $G_2 = G(j-1,n-k-j+1)$. Note that if $v \in V(H_k)$ in $G_2$, then $d_{G_2}(v) = d_{G_1}(v)$;  if $v \in V(K_j)$ in $G_2$, then $d_{G_2}(v) = d_{G_1}(v) - 1= j+k-2$;  if $v \in V(K_{n-k-j+1})$ in $G_2$, then $d_{G_2}(v) = d_{G_1}(v)+1 = n-j$.
By the definitions of $\prod_1$ and $\prod_2$, we have
\begin{eqnarray}
 \nonumber   \frac{\prod_1({G_1})}{\prod_1({G_2})}
&&=  \frac{\prod_{v \in V(K_j)}d(v)^2 \prod_{v \in V(H_k)}d(v)^2 \prod_{v \in V(K_{n-k-j})}d(v)^2}
{\prod_{v \in V(K_{j-1})}d(v)^2 \prod_{v \in V(H_k)}d(v)^2 \prod_{v \in V(K_{n-k-j+1})}d(v)^2}
 \\  \nonumber
&&= \frac{\big((j+k-1)^2\big)^j \big((n-j-1)^2 \big)^{n-k-j}}
{\big((j+k-2)^2\big)^{j-1} \big((n-j)^2 \big)^{n-k-j+1}}
  \\ \nonumber
&&= \bigg(\frac{ \frac{\big(j+(k-1)\big)^j}{\big((j-1) + (k-1)\big)^{j-1} }  }{ \frac{\big((n-j-k+1) + (k -1)\big)^{n-k-j+1}}{\big((n-j-k) + (k-1)\big)^{n-k-j} }}\bigg)^2.
 \end{eqnarray}

Since $2 \leq j \leq \frac{n-k}{2}$, we obtain $j \leq n-k-j < n-k-j+1$. By Proposition \ref{p1} and $k \geq 1$, we have $$\frac{\prod_1({G_1})}{\prod_1({G_2})} < 1,$$ that is, $\prod_1({G_1}) < \prod_1({G_2})$.

We can recursively use this process from $G_1$ to $G_2$, and obtain that
$$\prod_1(G(j,n-k-j)) < \prod_1(G(j-1,n-k-j+1)) < \prod_1(G(j-2,n-k-j+2)) < \cdots < \prod_1(G(1, n-k-1)).$$
Therefore, $\prod_1(G(j,n-k-j)) < \prod_1(G(1, n-k-1))$. Thus, we complete the proof. \qed
\end{proof}

\begin{lemma}\label{l4}
Let $G$ be a connected graph and $u, v \in V(G)$. Assume that $v_1, v_2, \ldots, v_s \in N(v)\setminus N(u)$,
$1 \leq s \leq d(v)$. Let $G' = G - \{vv_1, vv_2, \ldots, vv_s \} + \{uv_1, uv_2, \ldots, uv_s\}$.
If $d(u) \geq d(v)$ and $u$ is not adjacent to $v$, then
   $$\prod_2(G') > \prod_2(G). $$
\end{lemma}

\begin{proof}
By the concept of $\prod_2(G)$, we have
 \begin{eqnarray}
   \frac{\prod_2(G)}{\prod_2(G')} =
 \frac{d(u)^{d(u)} d(v)^{d(v)}}{(d(u)+s)^{d(u)+s}(d(v)-s)^{d(v)-s}}   = \frac{\big( \frac{d(u)^{d(u)}}{(d(u)+s)^{d(u)+s}}\big)}{\big( \frac{(d(v)-s)^{d(v)-s}}{d(v)^{d(v)}}\big)}.\nonumber
 \end{eqnarray}

By using $d(u) \geq d(v) > d(v) -s$ and   Proposition \ref{p2}, we obtain  $$\frac{\prod_2(G)}{\prod_2(G')} < 1,$$ which implies that
 $ \prod_2(G') > \prod_2(G)$. This shows the lemma. \qed
\end{proof}

\begin{lemma}\label{l5}
 If $k \geq 1$ and $2 \leq j \leq \frac{n-k}{2}$, we have $$\prod_2(G(j,n-k-j)) < \prod_2(G(1, n-k-1)).$$
\end{lemma}

\begin{proof}
Let $V(K_j) = \{v_1, v_2, \cdots, v_j\}$ and $ V(K_{n-k-j}) = \{u_1, u_2, \cdots, u_{n-k-j}\}$. Note that vertex set $\{ v_2, v_3, \cdots v_j \} \subset N(v_1) \cap V(K_j)$.
We define a new graph $G' = G(j,n-k-j) - \{v_1v_2, v_1v_3, \ldots, v_1v_j \} + \{u_1v_2, u_1v_3, \ldots, u_1v_j\}$. By $d(v_1) \leq d(u_1)$ and Lemma \ref{l4}, we have $\prod_2(G') \geq \prod_2(G(j,n-k-j))$.

Note that for $G'$, $v_1$ has neighbors in $V(H_k)$ only. Let $G'' = G' + \{v_iu_l, 2 \leq i \leq j, 1 \leq l \leq n-k-j$ and $v_iu_l \notin E(G')\}$.  By Lemma \ref{l2}, we have $\prod_2(G'') > \prod_2(G') \geq \prod_2(G)$.
Therefore $G'' \cong G(1, n-k-1)$ and $\prod_2(G(j,n-k-j)) < \prod_2(G(1, n-k-1))$. This completes the proof.  \qed
\end{proof}

\section{Extremal graphs with given connectivity}
In this section, the   maximal and minimal multiplicative Zagreb indices of graphs with connectivity at most $k$ in $\mathbb{V}_n^k$ and $\mathbb{E}_n^k$ are determined, and the corresponding extremal graphs have been characterized
in Theorems \ref{t1} and \ref{t2}.


\begin{theorem}\label{t1}
Let $G$ be a graph in $\mathbb{V}_n^k$. Then
\begin{eqnarray}&&\prod_1(G) \leq k^2 (n-k)^2k (n-2)^{2(n-k-1)}~ \text{and} \nonumber\\ \nonumber
&& \prod_2(G) \leq k^k (n-1)^{k(n-1)} (n-2)^{(n-2)(n-k-1)},
\end{eqnarray}
 where the equalities hold if and only if $G \cong K_n^k$.
\end{theorem}

\begin{proof}
Note that the degree sequence of $K_n^k$ is $k, \underbrace{n-2, n-2, \cdots, n-2}_{n-k-1}, \underbrace{n-1, n-1, \cdots, n-1}_k$. By the concepts of $\prod_1(G)$, $\prod_2(G)$ and routine calculations, we have
\begin{eqnarray}&&\prod_1(K_n^k) = k^2 (n-1)^{2k} (n-2)^{2(n-k-1)}~ \text{and} \nonumber\\ \nonumber
&& \prod_2(K_n^k) = k^k (n-1)^{k(n-1)} (n-2)^{(n-2)(n-k-1)}.
\end{eqnarray}
It suffices to prove that $\prod_1(G) \leq \prod_1(K_n^k)$ and $\prod_2(G) \leq \prod_2(K_n^k)$, and the equalities hold if and only if $G \cong K_n^k$.

If $k \geq n-1$, then $G \cong K_n^{n-1} \cong K_n$, and the theorem is true. If $1 \leq k \leq n-2$, then choose a graph $\overline{G}_1$ (resp. $\overline{G}_2$) in $\mathbb{V}_n^k$  such that $\prod_1(\overline{G}_1)$ (resp. $\prod_2(\overline{G}_2)$) is maximal.
Since $\overline{G}_i \ncong K_n$ with $i = 1, 2$, then $\overline{G}_i$ has a vertex cut set of size $k$. Let $V_i = \{v_{i1}, v_{i2}, \cdots, v_{ik}\}$ be the cut vertex set of $\overline{G}_i$.  Denoted $\omega(\overline{G}_i - V_i)$ by the number of  components of $\overline{G}_i - V_i$.  In order to prove our theorem, we start with several claims.

\begin{claim}\label{c1}
$\omega({\overline{G}_i - V_i}) = 2$ with $i =1, 2$.
\end{claim}
\noindent{\small \it Proof.}
We proceed to prove it by a contradiction. Assume that $\omega({\overline{G}_i - V_i}) \geq 3$ with $i = 1, 2$. Let $G_1, G_2, \cdots, G_{\omega({\overline{G}_i - V_i})}$ be the components of $\overline{G}_i - V_i$. Since $\omega({\overline{G}_i - V_i}) \geq 3$, then choose vertices $u \in V(G_1)$ and $v \in V(G_2)$. Then $V_i$ is still a $k$-vertex cut set of $\overline{G}_i+uv$. By Lemma \ref{l2}, we have $\prod_i(\overline{G}_i+uv) > \prod_i(\overline{G}_i)$, a contradiction to the choice of $\overline{G_i}$. Thus, this claim is proved.

Without loss of generality, suppose that $\overline{G}_i - V_i$ contains only two connected components, denoted by $G_{i1}$ and $G_{i2}$.

\begin{claim}\label{c2}
The induced graphs on $V(G_{i1}) \cup V_{i}$ and $V(G_{i2}) \cup V_{i}$  in $\overline{G}_i$ are complete subgraphs.
\end{claim}
\noindent{\small \it Proof.}
We use a contradiction to show it. Suppose that $\overline{G}_i[V(G_{i1}) \cup V_{i}]$  is not a complete subgraph of $\overline{G}_i$.
Then there exists an edge $uv \notin \overline{G}_i[V(G_{i1}) \cup V_{i}]$. Since $\overline{G}_i[V(G_{i1}) \cup V_{i}] + uv \in \mathbb{V}_{n,k}$, by Lemma \ref{l2}, we have $\prod_i(\overline{G}_i[V(G_{i1}) \cup V_{i}] +uv) > \prod_i(\overline{G}_i[V(G_{i1}) \cup V_{i}])$, which is a contadiction. This shows the claim.

By the above claims, we see that $G_{i1}$ and  $G_{i2}$ are complete subgraph of $\overline{G}_i$. Let $G_{i1} = K_{n'}$ and  $G_{i2}= K_{n''}$. Then we have $\overline{G}_i = K_{n'} \oplus  \overline{G}_i[V_i] \oplus K_{n''} $.

\begin{claim}\label{c3}
Either $n'=1$ or $n''=1$.
\end{claim}
\noindent{\small \it Proof.}
On the contrary, assume that $n', n'' \geq 2$. Without loss of generility, $n' \leq n''.$ For $\prod_i(G)$, by Lemmas \ref{l3} and \ref{l5}, we have a new graph $\overline{G'}_i = K_{1} \oplus  \overline{G}_i[V_i] \oplus K_{n-k-1} $ such that $\prod_i(\overline{G'}_i) > \prod_i(\overline{G}_i)$ and $\overline{G'}_i \in \mathbb{V}_n^k$. This is a contradition to the choice of $\overline{G}_i$.  Thus, either $n'=1$ or $n''=1$, and this claim is showed.

By Lemma \ref{l2}, $\prod_i(K_{1} \oplus  K_{|H_k|} \oplus K_{n-k-1}) > \prod_i( K_{1} \oplus  \overline{G}_i[V(H_k)] \oplus K_{n-k-1})$. Since  $\prod_i(K_n^k) = \prod_i(K_{1} \oplus  K_{|V_i|} \oplus K_{n-k-1})$, then  $\prod_i(K_n^k)$ is maximal and the theorem holds.
\end{proof}

Since $K_n^k \in \mathbb{E}_n^k \subset \mathbb{V}_n^k$,  the following result is immediate.
\begin{theorem}\label{t2}
Let $G$ be a graph in $\mathbb{E}_n^k$. Then
\begin{eqnarray}&&\prod_1(G) \leq k^2 (n-k)^2k (n-2)^{2(n-k-1)}  \text{and} \nonumber\\ \nonumber
&& \prod_2(G) \leq k^k (n-1)^{k(n-1)} (n-2)^{(n-2)(n-k-1)},
\end{eqnarray}
 where the equalities hold if and only if $G \cong K_n^k$.
\end{theorem}

In the rest of this Section, we consider the minimal mutiplicative Zagreb indices of graphs $G$ in $\mathbb{V}_n^k$ and $\mathbb{E}_n^k$. By Proposition \ref{p4} (ii), $G$ is a tree with $n$ vertices.  By Lemma \ref{l1} and routine calculations, we have

\begin{theorem} \label{t3}
Let $G$ be a graph in $\mathbb{V}_n^k$. Then
$$\prod_1(G) \geq (n-1)^2~\text{and} ~ \prod_2(G) \geq 4^{n-2},$$ where the equalities hold if and only if $G \cong S_n$ and $G \cong P_n$, respectively.
\end{theorem}

Note that  $P_n, S_n \in \mathbb{E}_n^k \subset \mathbb{V}_n^k$, then the following theorem is obvious.

\begin{theorem} \label{t4}
Let $G$ be a graph in $\mathbb{E}_n^k$. Then
$$\prod_1(G) \geq (n-1)^2~\text{and} ~ \prod_2(G) \geq 4^{n-2},$$ where the equalities hold if and only if $G \cong S_n$ and $G \cong P_n$, respectively.
\end{theorem}

\section{Extremal graphs with given number of pendant vertices}
Let $\mathbb{G}_n^p$ be the set of graphs with $p \geq 2$ pendant vertices. In this section, the maximal and minimal multiplicative Zagreb indices of graphs with $p$ pendant vertices in $\mathbb{G}_n^p$ are determined, and the corresponding extremal graphs shall be characterized
in Theorems \ref{p11} and \ref{p12}.

Before exhibiting the main results of the section, we list some notations which will be used in the sequel.
Clearly, if $G\in \mathbb{G}_n^p$, then there be a connected subgraph $H_1$ with order $n-p$ for which $G$ can reconstructed by linking $p$ vertices to some vertices $H_1$. Especially, since $H_1$ is connected, it has two extremal cases, i.e., $H_1\cong K_{n-p}$ and $H_1\cong T_{n-p}.$ Let $\mathcal{A}^1_n$ and $\mathcal{A}^2_n$ be the
two graph sets such that its element with the sequence $(p,\underbrace{2,\ldots,2}_{n-p-1},\underbrace{1,\ldots,1}_{p})$ and $(\underbrace{k+1,\ldots,k+1}_{r},\underbrace{k,\ldots,k}_{n-p-r},
\underbrace{1,\ldots,1}_{p}),$ where  $2n-p-2=k(n-p)+r$ with $k\geq 2$ and $0\leq r\leq n-p-1.$
 Let $T$ be a tree, and $v\in V(T)$ with $d(v)=k$. Note that $T-v$ has $k$ components, for each component associated with $v$, we call it as a \emph{branch} of $v.$ We notice that graph $G_a$ meets $|n_i-n_j|\leq 1$ for $1\leq i,j\leq n-p$, see Fig.2. Since $\sum_i^{n-p}n_i=p$.  There are two integers $\ell$ and $t$ for which $p=\ell(n-p)+t$ with $\ell\geq 0$ and $0\leq t\leq n-p-1.$ In other words, $G_a$ has the sequence $(\underbrace{\ell+1,\ldots,\ell+1}_{t},
 \underbrace{\ell,\ldots,\ell}_{n-p-t},
 \underbrace{1,\ldots,1}_{p}).$

\vspace{5mm}
\begin{center}
\begin{picture}(217,50)\linethickness{0.8pt}
\cbezier(3,18)(5.4,-5.9)(73.6,-5.9)(76,18)
\cbezier(3,18)(5.4,41.9)(73.6,41.9)(76,18)
\Line(0,49.5)(14,25)
\Line(16,49.5)(14,25)
\Line(24.5,50)(30.5,24.5)
\Line(40,49.5)(30.5,24.5)
\Line(52.5,48)(58.5,23.5)
\Line(70.5,47.5)(58.5,23.5)
\cbezier(145,17.3)(147.3,-5.6)(214.7,-5.6)(217,17.3)
\cbezier(145,17.3)(147.3,40.1)(214.7,40.1)(217,17.3)
\Line(158.5,46)(172,21.5)
\Line(175.5,46)(172,21.5)
\Line(198,46)(172,21.5)
\put(166.1,6.6){\footnotesize$K_{n-p}$}
\put(20.6,4.6){\footnotesize$K_{n-p}$}
\put(27.2,15.1){\footnotesize$v_2$}
\put(51.2,13.6){\footnotesize$v_{n-p}$}
\put(6.2,16.6){\footnotesize$v_1$}
\put(80,0){\footnotesize$G_{a}$}
\put(140,0){\footnotesize$G_{s}$}
\put(14,25){\circle*{4}}
\put(30.5,24.5){\circle*{4}}
\put(38,24.5){\circle*{2}}
\put(43,24.5){\circle*{2}}
\put(48,24.5){\circle*{2}}
\put(58.5,23.5){\circle*{4}}
\put(0,49.5){\circle*{4}}
\put(16,49.5){\circle*{4}}
\put(5.5,48){\circle*{2}}
\put(13,48){\circle*{2}}
\put(9.5,48){\circle*{2}}
\put(24.5,50){\circle*{4}}
\put(40,49.5){\circle*{4}}
\put(28,48){\circle*{2}}
\put(31,48){\circle*{2}}
\put(35,48){\circle*{2}}
\put(52.5,48){\circle*{4}}
\put(70.5,47.5){\circle*{4}}
\put(57,48){\circle*{2}}
\put(61,48){\circle*{2}}
\put(65,48){\circle*{2}}
\put(158.5,46){\circle*{4}}
\put(172,21.5){\circle*{4}}
\put(175.5,46){\circle*{4}}
\put(182,46){\circle*{2}}
\put(187,46){\circle*{2}}
\put(192,46){\circle*{2}}
\put(198,46){\circle*{4}}
\put(8,55){\makebox(0,0){$\overbrace{\rule{5mm}{0mm}}^{\makebox(0,0){\footnotesize$n_1$}}$}}
\put(32,55){\makebox(0,0){$\overbrace{\rule{5mm}{0mm}}^{\makebox(0,0){\footnotesize$n_2$}}$}}
\put(61,55){\makebox(0,0){$\overbrace{\rule{5mm}{0mm}}^{\makebox(0,0){\footnotesize$n_{n-p}$}}$}}
\put(177,55){\makebox(0,0){$\overbrace{\rule{15mm}{0mm}}^{\makebox(0,0){\footnotesize$p$}}$}}
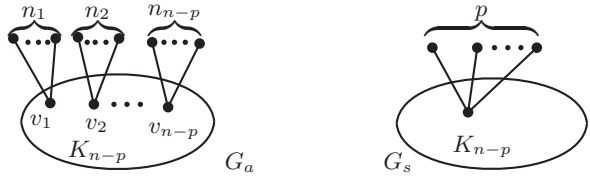
\put(107,-20){\makebox(0,0){Figure.2 The graphs $G_a$ with $|n_i-n_j|\leq 1$ for $1\leq i,j\leq n-p$ and $G_s$.}}
\end{picture}
\end{center}
\vspace{5mm}
\begin{theorem} \label{p11}
Let $G$ be a graph in $\mathbb{G}_n^p$. Then
$$\prod_1(G)\leq (n+\ell-p)^{2t}(n+\ell-p-1)^{2(n-p-\ell)}~\text{and}~\prod_2(G)\leq (n-1)^{n-1}(n-p-1)^{(n-p-1)^2},$$
where the equalities hold if and only if $G \cong G_{a}$ and $G\cong G_{s}$, respectively(see, Fig. 2).
\end{theorem}
\proof Suppose that $G\in \mathbb{G}_n^p$ such that $G$ has the maximum value with respect to $\prod_1$ and $\prod_2$. According to properties of $\prod_i$ for $i=1,2$, if $G+e\in \mathbb{G}_n^p$, we obtain that $\prod_i(G)<\prod_i(G+e).$ Hence the subgraph $H_1$ of $G$ with order $n-p$ is the complete graph $K_{n-p}.$ Labeling the vertices of $H_1$ as $v_1,v_2,\ldots,v_{n-p}$, let $n_i$ be the number of pendant vertices who link with $v_i$, for $i=1,2,\ldots,n-p$. We firstly show the upper bound of $\prod_1.$

Assume that there are two vertices $v_i$ and $v_j$ of $H_1$ such that $|n_i-n_j|\geq 2.$
With loss of generality, set $n_i-n_j\geq 2.$ Let $G'$ be the new graph from $G$ by deleting one pendent vertex and adding it to $v_j.$ Note that $d_{G'}(v_i)=d_{G}(v_i)-1$ and $d_{G'}(v_j)=d_{G}(v_j)+1.$  For convenience, we write $d$ instead of $d_{G}$.  Observe that
\begin{equation*}
\begin{split}
\frac{\prod_1(G')}{\prod_1(G)}=\frac{d^2_{G'}(v_i)d^2_{G'}(v_j)}{d^2_{G}(v_i)d^2_{G}(v_j)}
&=\frac{(d(v_i)-1)^2(d(v_j)+1)^2}{d^2(v_i)d^2(v_j)}\\
&=1+\frac{2d(v_i)d(v_j)(d(v_i)-d(v_j)-2)+(d(v_i)-1)^2+d^2(v_j)+2d(v_j)}{d^2(v_i)d^2(v_j)}\\
&\geq 1+\frac{d^2(v_j)+2d(v_j)}{d^2(v_i)d^2(v_j)}>1.
\end{split}
\end{equation*} Consequently, $\prod_1(G')>\prod_1(G)$ which contradicts  with the choice of $G$.
Hence, for any pair $v_i$ and $v_j$ of $H_1$, $|n_i-n_j|\leq 1.$ In other words, $G\cong G_a.$  Clearly, by routine calculation, $\prod_1(G_a)= (n+\ell-p)^{2t}(n+\ell-p-1)^{2(n-p-\ell)}.$

We now verify the
upper bound of $\prod_2.$ In order to obtain the maximum of $\prod_2$, it is sufficient to
 show the following claim.
\begin{claim}
 The number of the vertices in $H_1$ who possesses pendent vertex  is one.
\end{claim}
\proof
 Assume that $G$ has at least two vertices, such as, $v_i$ and $v_j(\,\text{with}\, d(v_i)\geq d(v_j))$, which have pendents.
 Denote by $G''$  the graph obtained from $G$ by deleting one pendent of $v_j$ and adding to $v_i.$

 So $d_{G''}(v_i)-d_{G''}(v_j)\geq 2.$

 According to Proposition \ref{p1}, we observe that
 \begin{equation*}
 \begin{split}
 \frac{\prod_2(G'')}{\prod_2(G)}=\frac{d_{G''}(v_i)^{d_{G''}(v_i)}d_{G''}(v_j)^{d_{G''}(v_j)}}
 {d_{G}(v_i)^{d_{G}(v_i)}d_{G}(v_j)^{d_{G}(v_j)}}
 &=\frac{(d(v_i)+1)^{d(v_i)+1}(d(v_j)-1)^{d(v_j)-1}}
 {d(v_i)^{d(v_i)}d(v_j)^{d(v_j)}}\\
 &>1(\text{by}\ \text{setting}\  d(v_i)\triangleq d(v_j)-1).
 \end{split}
 \end{equation*}So, $\prod_2(G'')>\prod_2(G)$, a contradiction. Therefore, the claim is holds.\qed

 Clearly, $G$ is the graph such that $\prod_1$ has maximum if and only if $G\cong G_a$, and $G$ is maximal graph regarding $\prod_2$ if and only if $G\cong G_s.$
 By direct calculation, we have $\prod_1(G_a)=(n+\ell-p)^{2t}(n+\ell-p-1)^{2(n-p-\ell)}$ and  $\prod_1(G_s)=(n-1)^{n-1}(n-p-1)^{(n-p-1)^2}.$

 Therefore, we complete the proof. \qed

\vspace{5mm}
\begin{center}
\begin{picture}(277.5,55.9)\linethickness{0.8pt}
\put(51.3,8){$G$}
\put(207.8,8){$G_1$}
\put(193,19){$v_1$}
\put(234.2,19){$v_2$}
\put(30.2,19){$v_1$}
\put(77.2,19){$v_2$}
\Line(104,27.4)(123.5,27.4)
\Line(83.5,27.4)(104,27.4)
\Line(95.5,43.4)(83.5,27.4)
\Line(86.5,46.4)(83.5,27.4)
\Line(71,45.9)(83.5,27.4)
\Line(37,27.4)(83.5,27.4)
\Line(42.5,46.9)(37,27.4)
\Line(51.5,46.4)(37,27.4)
\Line(27,46.9)(37,27.4)
\Line(18.5,27.9)(37,27.4)
\Line(196,46.4)(191,26.9)
\Line(181,46.4)(191,26.9)
\Line(258,26.9)(277.5,26.9)
\Line(237.5,26.9)(258,26.9)
\Line(191,26.9)(237.5,26.9)
\Line(170.5,26.9)(191,26.9)
\Line(152.5,26.9)(170.5,26.9)
\Line(0,27.9)(18.5,27.9)
\cbezier(95,46.4)(95.2,40.4)(101.8,40.4)(102,46.4)
\cbezier(95,46.4)(95.2,52.4)(101.8,52.4)(102,46.4)
\cbezier(64.5,48.9)(64.7,42.9)(71.8,42.9)(72,48.9)
\cbezier(64.5,48.9)(64.7,54.9)(71.8,54.9)(72,48.9)
\cbezier(174,49.9)(174.3,43.9)(181.7,43.9)(182,49.9)
\cbezier(174,49.9)(174.3,55.9)(181.7,55.9)(182,49.9)
\cbezier(84,50.9)(84.2,44.9)(90.3,44.9)(90.5,50.9)
\cbezier(84,50.9)(84.2,56.9)(90.3,56.9)(90.5,50.9)
\cbezier(50.5,50.7)(50.7,45)(57.3,45)(57.5,50.7)
\cbezier(50.5,50.7)(50.7,56.3)(57.3,56.3)(57.5,50.7)
\cbezier(39.5,51.4)(39.7,45.4)(46.3,45.4)(46.5,51.4)
\cbezier(39.5,51.4)(39.7,57.4)(46.3,57.4)(46.5,51.4)
\cbezier(20,50.4)(20.3,44.4)(27.7,44.4)(28,50.4)
\cbezier(20,50.4)(20.3,56.4)(27.7,56.4)(28,50.4)
\cbezier(195,50.4)(195.2,44.4)(201.8,44.4)(202,50.4)
\cbezier(195,50.4)(195.2,56.4)(201.8,56.4)(202,50.4)
\Line(191,26.9)(206.5,43.4)
\cbezier(205.5,46.9)(205.7,41.6)(212.3,41.6)(212.5,46.9)
\cbezier(205.5,46.9)(205.7,52.2)(212.3,52.2)(212.5,46.9)
\Line(191,26.9)(177,15.9)
\cbezier(175.5,15.4)(175.5,15.4)(175.5,15.4)(175.5,15.4)
\cbezier(175.5,15.4)(175.5,15.4)(175.5,15.4)(175.5,15.4)
\cbezier(170.5,15.7)(170.7,10)(177.3,10)(177.5,15.7)
\cbezier(170.5,15.7)(170.7,21.3)(177.3,21.3)(177.5,15.7)
\Line(191,26.9)(192,14.4)
\cbezier(186.5,10.4)(186.7,4.4)(193.8,4.4)(194,10.4)
\cbezier(186.5,10.4)(186.7,16.4)(193.8,16.4)(194,10.4)
\put(123.5,27.4){\circle*{4}}
\put(104,27.4){\circle*{4}}
\put(60.5,27.4){\circle*{4}}
\put(83.5,27.4){\circle*{4}}
\put(37,27.4){\circle*{4}}
\put(214.5,26.9){\circle*{4}}
\put(277.5,26.9){\circle*{4}}
\put(258,26.9){\circle*{4}}
\put(237.5,26.9){\circle*{4}}
\put(191,26.9){\circle*{4}}
\put(170.5,26.9){\circle*{4}}
\put(152.5,26.9){\circle*{4}}
\put(237.5,27.9){\circle*{4}}
\put(95.5,43.4){\circle*{4}}
\put(86.5,46.4){\circle*{4}}
\put(71,45.9){\circle*{4}}
\put(42.5,46.9){\circle*{4}}
\put(51.5,46.4){\circle*{4}}
\put(27,46.9){\circle*{4}}
\put(18.5,27.9){\circle*{4}}

\put(196,46.4){\circle*{4}}
\put(181,46.4){\circle*{4}}
\put(0,27.9){\circle*{4}}

\put(190.5,45){\circle*{2}}
\put(187.5,45){\circle*{2}}
\put(184.5,45){\circle*{2}}

\put(82,44.9){\circle*{2}}
\put(79,44.9){\circle*{2}}
\put(76,44.9){\circle*{2}}

\put(37.5,45.9){\circle*{2}}
\put(34,45.9){\circle*{2}}
\put(30.5,45.9){\circle*{2}}

\put(206.5,43.4){\circle*{4}}
\put(177,15.9){\circle*{4}}
\put(192,14.4){\circle*{4}}

\put(182,16.9){\circle*{2}}
\put(185,16.9){\circle*{2}}
\put(188,16.9){\circle*{2}}
\put(40,60){\makebox(0,0){$\overbrace{\rule{10mm}{0mm}}^{\makebox(0,0){\footnotesize$n_1$}}$}}
\put(85,60){\makebox(0,0){$\overbrace{\rule{12mm}{0mm}}^{\makebox(0,0){\footnotesize$n_2$}}$}}
\put(195,60){\makebox(0,0){$\overbrace{\rule{12mm}{0mm}}^{\makebox(0,0){\footnotesize$n_1$}}$}}
\put(182,6){\makebox(0,0){$\underbrace{\rule{6mm}{0mm}}_{\makebox(0,0){\footnotesize$n_2$}}$}}
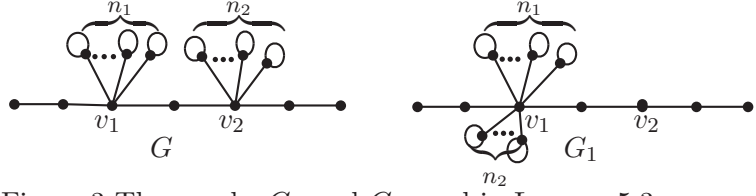
\put(120,-10){\makebox(0,0){Figure.3 The graphs $G_1$ and $G_2$ used in Lemma \ref{51}.}}
\end{picture}
\end{center}

\vspace{5mm}
\begin{center}
\begin{picture}(277.5,55.9)\linethickness{0.8pt}
\Line(152.5,26.9)(170.5,26.9)
\Line(170.5,26.9)(191,26.9)
\Line(191,26.9)(237.5,26.9)
\Line(237.5,26.9)(258,26.9)
\Line(258,26.9)(277.5,26.9)
\Line(181,46.4)(191,26.9)
\Line(196,46.4)(191,26.9)
\Line(225,45.4)(237.5,26.9)
\Line(239,46.4)(237.5,26.9)
\Line(247,45.4)(237.5,26.9)
\cbezier(174,49.9)(174.3,43.9)(181.7,43.9)(182,49.9)
\cbezier(174,49.9)(174.3,55.9)(181.7,55.9)(182,49.9)
\cbezier(195,50.4)(195.2,44.4)(201.8,44.4)(202,50.4)
\cbezier(195,50.4)(195.2,56.4)(201.8,56.4)(202,50.4)
\cbezier(218.5,48.4)(218.7,42.4)(225.8,42.4)(226,48.4)
\cbezier(218.5,48.4)(218.7,54.4)(225.8,54.4)(226,48.4)
\cbezier(235.5,50.9)(235.7,44.9)(241.8,44.9)(242,50.9)
\cbezier(235.5,50.9)(235.7,56.9)(241.8,56.9)(242,50.9)
\cbezier(246.5,49.4)(246.7,43.4)(253.3,43.4)(253.5,49.4)
\cbezier(246.5,49.4)(246.7,55.4)(253.3,55.4)(253.5,49.4)
\Line(0,27.9)(18.5,27.9)
\cbezier(20,50.4)(20.3,44.4)(27.7,44.4)(28,50.4)
\cbezier(20,50.4)(20.3,56.4)(27.7,56.4)(28,50.4)
\Line(18.5,27.9)(37,27.4)
\Line(27,46.9)(37,27.4)
\cbezier(39.5,51.4)(39.7,45.4)(46.3,45.4)(46.5,51.4)
\cbezier(39.5,51.4)(39.7,57.4)(46.3,57.4)(46.5,51.4)
\Line(51.5,46.4)(37,27.4)
\Line(42.5,46.9)(37,27.4)
\Line(37,27.4)(83.5,27.4)
\cbezier(50.5,50.7)(50.7,45)(57.3,45)(57.5,50.7)
\cbezier(50.5,50.7)(50.7,56.3)(57.3,56.3)(57.5,50.7)
\cbezier(64.5,48.9)(64.7,42.9)(71.8,42.9)(72,48.9)
\cbezier(64.5,48.9)(64.7,54.9)(71.8,54.9)(72,48.9)
\Line(71,45.9)(83.5,27.4)
\cbezier(84,50.9)(84.2,44.9)(90.3,44.9)(90.5,50.9)
\cbezier(84,50.9)(84.2,56.9)(90.3,56.9)(90.5,50.9)
\Line(86.5,46.4)(83.5,27.4)
\Line(95.5,43.4)(83.5,27.4)
\Line(83.5,27.4)(104,27.4)
\cbezier(95,46.4)(95.2,40.4)(101.8,40.4)(102,46.4)
\cbezier(95,46.4)(95.2,52.4)(101.8,52.4)(102,46.4)
\Line(104,27.4)(123.5,27.4)
\Line(259.5,41.4)(237.5,27.9)
\cbezier(259.5,44.2)(259.8,39.2)(267.7,39.2)(268,44.2)
\cbezier(259.5,44.2)(259.8,49.1)(267.7,49.1)(268,44.2)
\put(234.2,19){$v_1$}
\put(187.7,19){$v_2$}
\put(30.2,19){$v_2$}
\put(77.2,19){$v_1$}
\put(207.8,8){$G_2$}
\put(51.3,8){$G$}
\put(152.5,26.9){\circle*{4}}
\put(170.5,26.9){\circle*{4}}
\put(191,26.9){\circle*{4}}
\put(237.5,26.9){\circle*{4}}
\put(258,26.9){\circle*{4}}
\put(277.5,26.9){\circle*{4}}
\put(214.5,26.9){\circle*{4}}
\put(181,46.4){\circle*{4}}
\put(196,46.4){\circle*{4}}
\put(225,45.4){\circle*{4}}
\put(239,46.4){\circle*{4}}
\put(247,45.4){\circle*{4}}
\put(0,27.9){\circle*{4}}
\put(18.5,27.9){\circle*{4}}
\put(37,27.4){\circle*{4}}
\put(27,46.9){\circle*{4}}
\put(51.5,46.4){\circle*{4}}
\put(42.5,46.9){\circle*{4}}
\put(83.5,27.4){\circle*{4}}
\put(60.5,27.4){\circle*{4}}
\put(71,45.9){\circle*{4}}
\put(86.5,46.4){\circle*{4}}
\put(95.5,43.4){\circle*{4}}
\put(104,27.4){\circle*{4}}
\put(123.5,27.4){\circle*{4}}
\put(259.5,41.4){\circle*{4}}
\put(237.5,27.9){\circle*{4}}

\put(30.5,45.9){\circle*{2}}
\put(34.5,45.9){\circle*{2}}
\put(38.5,45.9){\circle*{2}}

\put(76,45){\circle*{2}}
\put(79.5,45){\circle*{2}}
\put(83,45){\circle*{2}}

\put(185,45){\circle*{2}}
\put(188,45){\circle*{2}}
\put(191,45){\circle*{2}}

\put(229,45){\circle*{2}}
\put(232.5,45){\circle*{2}}
\put(236,45){\circle*{2}}
\put(40,60){\makebox(0,0){$\overbrace{\rule{10mm}{0mm}}^{\makebox(0,0){\footnotesize$n_2$}}$}}
\put(85,60){\makebox(0,0){$\overbrace{\rule{12mm}{0mm}}^{\makebox(0,0){\footnotesize$n_1$}}$}}
\put(189,60){\makebox(0,0){$\overbrace{\rule{8mm}{0mm}}^{\makebox(0,0){\footnotesize$n_2-1$}}$}}
\put(243,60){\makebox(0,0){$\overbrace{\rule{14mm}{0mm}}^{\makebox(0,0){\footnotesize$n_1+1$}}$}}
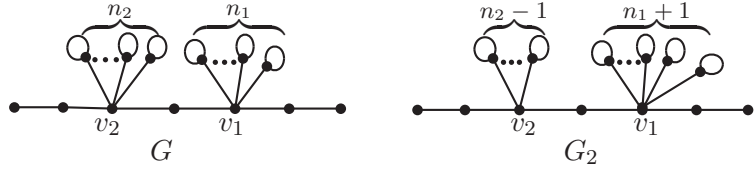
\put(140,-10){\makebox(0,0){Figure.4 The graphs $G_3$ and $G_4$ used in Lemma \ref{52}}}
\end{picture}
\end{center}

\begin{lemma}\label{51}
  If $G$ and $G_1$ are two graphs as shown in Fig. 3, and $G_1$ is regarded as the graph obtained from $G$ by transferring $n_2$ branches of $v_2$ to $v_1.$ Then $\prod_1(G_1)<\prod_1(G).$
\end{lemma}
\proof  We always suppose $d_G(v_1)\geq d_G(v_2)$ (If not, exchanging the signs of $v_1$ and $v_2$). Clearly, $d_G(v_1)=n_1+2$ and $d_G(v_2)=n_2+2.$ In view of the property of $\prod_1$ and Proposition \ref{p3}, we have
\begin{equation*}
\begin{split}
\frac{\prod_1(G_1)}{\prod_1(G)}
&=\frac{d_{G_1}(v_1)^2d_{G_1}(v_2)^2}{d_{G}(v_1)^2d_G(v_2)^2}\\
&=\frac{(n_1+n_2+2)^22^2}{(n_1+2)^2(n_2+2)^2}\\
&<1.
\end{split}
\end{equation*} Hence, the proof is complete.\qed
\begin{lemma}\label{52}
Let $G$ and $G_2$ be two graphs as shown in Fig. 4, and $G_2$ is considered as the graph obtained from $G$ by deleting one branch of $v_2$ and adding to $v_1$. If $d_G(v_2)-d_G(v_1)\geq 2$. Then $\prod_1(G_2)<\prod_1(G).$
\end{lemma}
\proof Note that $d_G(v_1)=n_1+2$ and $d_G(v_2)=n_2+2.$ Since $d_G(v_2)-d_G(v_1)\geq 2$. According to the Proposition \ref{p2} and the property of $\prod_2$, we find that
\begin{equation*}
\begin{split}
\frac{\prod_2(G_2)}{\prod_2(G)}
&=\frac{d_{G_2}(v_1)^{d_{G_2}(v_1)}d_{G_2}(v_2)^{d_{G_2}(v_2)}}
{d_{G}(v_1)^{d_{G}(v_1)}d_{G}(v_2)^{d_{G}(v_2)}}\\
&=\frac{(n_1+3)^{n_1+3}(n_2+1)^{n_2+1}}{(n_1+2)^{n_1+2}(n_2+2)^{n_2+2}}\\
&<\frac{(n_1+3)^{n_1+3}(n_1+1+1)^{n_1+1+1}}{(n_1+2)^{n_1+2}(n_1+1+2)^{n_1+1+2}}=1.
\end{split}
\end{equation*} Therefore, we finish the proof. \qed

\begin{theorem} \label{p12}
Let $G$ be a graph in $\mathbb{G}_n^p$. Then
$$\prod_1(G)\geq p^22^{2(n-p-1)}~\text{and}~\prod_2(G)\geq (k+1)^{r(k+1)}k^{k(n-p-r)},$$
where the equalities hold if and only if $G \in \mathcal {A}^1_n$ and $G \in \mathcal{A}^2_n$, respectively.
\end{theorem}

\proof   Let $G^{*}$ be the minimal graph with respect to $\prod_1$ and $\prod_2$,  respectively. Obviously, $H_1$ of $G^{*}$ is a tree $T_{n-p}$ through $\prod_i(G-e)<\prod_i(G)$ for $i=1,2.$
 We first consider the lower bound of $\prod_1.$

If $G^{*}$ just has one vertex whose degree is more than three. According to a property of tree, $G^{*}$ has $p$ pendent vertices. Then $G^{*}\in \mathcal{A}_n$. Otherwise, assume that $G^{*}$ has at least two vertices with degree more than two (they belong to $H_1$.), such as $v_i$ and $v_j$ (suppose $d_{G^*}(v_i)\geq d_{G^*}(v_j)$).  Since two vertices of a tree have a unique path through them. Let $P_t$ be a maximal path via $v_i$ and $v_j$. We call graph $G_1$ be the new graph obtained from $G^*$ by
deleting $d_{G^*}(v_j)-2$ branches of $v_j$ and linking to $v_i.$ In terms of Lemma \ref{51}, it is easy to deduce that $\prod_1(G_1)<\prod_1(G^*).$ A contradiction finish the proof of the part.

Next, we discuss the lower bound of $\prod_2$. Since $G^*$ has $p$ pendents. Labeling all vertices of $H_1$ as $v_1,v_2,\ldots,v_{n-p},$ $d_{G^*}(v_i)\geq 2$ for $1\leq i\leq n-p.$
We claim that $|d_{G^*}(v_i)-d_{G^*}(v_j)|\leq 1$ for all $1\leq i,j\leq n-p.$ If not, there exists at least two vertices in $H_1$, such as, $v_{i_0}$ and $v_{j_0}$, such that $d_{G^*}(v_{i_0})-d_{G^*}(v_{j_0})\geq 2.$ Considering the new graph $G_2$ obtained from $G^*$ by transferring one branch of $v_{i_0}$ to $v_{j_0},$ by means of Lemma \ref{52}, we get $\prod_2(G_2)<\prod_2(G^*),$ which is contradicted with the choice of $G^*.$ Observe that $\sum_{v\in V(G*)}(d_{G*}(v))=2(n-1)$, and $\sum_{v\in V(H_1)}(d_{G*}(v))=2(n-1)-p.$
Hence, there exist two integers $k,r$ such that $2(n-1)-p=k(n-p)+r$, where $k\geq 2$ and $0\leq r\leq n-p-1.$

Combining the above discussion, we deduce that $G^*$ belongs to $\mathcal {A}^1_n$ for $\prod_1,$ and $G^*$ belongs to $\mathcal {A}^2_n$ with respect to $\prod_2.$
 From the definition of $\mathcal {A}^1_n$ and $\mathcal {A}^2_n$, through a direct calculation, we have $\prod_1(G^*)= p^22^{2(n-p-1)}$ and $\prod_2(G^*)=(k+1)^{r(k+1)}k^{k(n-p-r)}$.

 Therefore, we finish the proof. \qed

\vskip4mm\noindent{\bf\large Acknowledgments}

\noindent This work was supported by National Natural Science Foundation of China (Grant Nos.11401348 and 11561032), Shandong Provincial Natural Science Foundation of China(No. ZR2019MA012),and supported by Postdoctoral Science Foundation of China.

%
%
%
%
%
%
%
%
%

\end{document}